\documentclass[12pt]{amsart}

\usepackage{amsfonts}

\usepackage{mathtext}
\usepackage[cp1251]{inputenc}
\usepackage[T2A]{fontenc}
\RequirePackage{srcltx}
\usepackage[english]{babel}

\usepackage[dvips]{graphicx}
\usepackage{amsmath}
\usepackage{amssymb}
\usepackage{amsxtra}
\usepackage{latexsym}
\usepackage{ifthen}
\usepackage{mathrsfs}

\theoremstyle{plain}
\newtheorem{theorem}{Теорема}[section]
\newtheorem{lemma}{Lemma}[section]

\newtheorem{corollary}{Corollary}[section]

\theoremstyle{definition}

\def\){\right)}
\def\({\left(}

\def\N{{{\Bbb N}}}
\def\Z{{{\Bbb Z}}}
\def\T{{{\Bbb T}}}

\begin{document}

\title[Almost everywhere summability of Fourier series]{Almost everywhere summability of Fourier series with indicating
the set of convergence}

\author{R. M. Trigub}

\email{roald.trigub@gmail.com}



\date{}

\maketitle

\begin{abstract}
The following problem is studied in this paper: Which multipliers $\{\lambda_{k, n}\}$ ensure
the convergence, as $n\to \infty$, of the linear means of the Fourier series of functions $f\in L_1[-\pi, \pi]$
$$
\sum_{k=-\infty}^\infty \lambda_{k, n}\hat{f}_k e^{ikx},
$$
where $\widehat{f}_k$ is the $k$-th Fourier coefficient, at a point at which the derivative of the function $\int_0^x f$ exists.
A criterion for the convergence of the $(C, 1)$-means ($\lambda_{k, n}=(1-\frac {|k|}{n+1})_+$) is found, while in the general
case $\lambda_{k, n}=\phi(\frac {k}{n+1})$ a sufficient condition is derived for the convergence at all such points (that is,
almost everywhere). The answer is given in terms of the belonging of $\phi(x)$ and $x\phi'(x)$ to the Wiener algebra of absolutely
convergent Fourier integrals. The obtained results are supplemented by some examples.
\end{abstract}

\vskip4mm

\medskip

\section{Введение}
\medskip

Ряд Фурье $2\pi$-периодической функции $f\in L_1(\mathbb T),  \T=[-\pi, \pi]$,  будем
записывать в виде
$$
f\sim\sum_{k=-\infty}^\infty \widehat{f}_k e_k, \quad e_k=e^{ikx},  \quad \widehat
{f}_k=\frac{1}{2\pi}\int_{-\pi}^\pi f(t)e^{-ikt}dt.
$$
Как известно, он может расходиться всюду (А.Н. Колмогоров [1]), тогда как ряд Фурье любой
функции $f\in L_p(\mathbb T)$, $p\in (1, +\infty),$ сходится почти всюду (Л. Карлесон, Р. Хант;
см.  [2]).  Сравнительно недавно получены существенные усиления этих
 результатов в их сближении (см.  [3],  [4]).

 С другой стороны, уже давно изучают сходимость при $n\to \infty$ линейных средних рядов
 Фурье вида
 $$
 \Lambda_n(f;x)\sim\sum_{k=-\infty}^\infty \lambda_{k, n}\widehat{f}_k e^{ikx}                   \eqno (1)
 $$
 в зависимости от множителей  $\lambda_{k, n}$ (см. напр.,  [5] т. I, гл. III, [6]).  Это свёртки
функции $f$ с ядрами $K_n$:
$$
\Lambda_{n}(f;x)=\frac{1}{2\pi}\int_{-\pi}^{\pi}f(x-t)K_n(t)dt, \quad K_n(t)\sim
\sum_{-\infty}^\infty \lambda_{k, n}e^{ikt}.
$$
Здесь важны вопросы о сходимости к $f$ по норме $L_1(\mathbb T)$ и о поточечной
сходимости почти всюду.

В случае $\lambda_{k, n}=\phi(\frac{k}{n+1})$,  где $\phi$ -- ограниченная и непрерывная почти
всюду на $\mathbb R$ функция, и сходимости по норме $L_1(\mathbb T)$ (или $C(\mathbb
T)$) имеется следующий критерий,  т.е.,  необходимое и достаточное условие одновременно:
функция $\phi$ (после исправления по непрерывности) является преобразованием Фурье
конечной на $R$ комплекснозначной борелевской меры и $\phi(0)=1$ (см.  [7],  8.1.2).

А. Лебег ввёл точки $x$ ($l$-точки),  для которых существует $l_f(x)$ с условием
$$
\lim_{|h|\to 0}\frac1h\int_0^h|f(x+t)-l_f(x)|dt=0,
$$
и доказал, что для любой функции $f\in L_1(\mathbb T)$ почти все точки являются её
точками Лебега.  Кроме того, он доказал, что средние арифметические частных сумм Фурье
($(C, 1)$-средние)
$$ \sigma_n(f)=\frac{1}{n+1}\sum_{m=0}^{n} s_m(f)=\sum_{k=-n}^n
\(1-\frac{|k|}{n+1}\)\widehat{f}_k e_k,  \eqno (2)
$$
$$
s_m(f)=\sum_{k=-m}^m \widehat{f}_ke_k
$$
сходятся к $l_f(x)$ во всех $l$-точках любой функции $f\in L_1(\mathbb T)$ (см.
$[5],  [6],  [8]$).  Для сходимости во всех точках Лебега имеется критерий в
терминах "горбатой" мажоранты модуля ядра $K_n$ [9],  а в случае
$\lambda_{k, n}=\phi(\frac{k}{n+1})$ - в терминах преобразования Фурье $\phi$
(точнее,  $\phi(0)=1$ и $\phi$ принадлежит алгебре $A^*$,  определение которой отличается от
определения $A(\mathbb R)$ (см. ниже (4)) тем,  что вместо $g\in L_1$ должно быть
${\rm ess}\, \sup_{|x|\ge t}|g(x)|\in L_1[0, +\infty)$  ([7],  8.1.3).

В настоящей статье будем изучать сходимость во всех точках $x$,  в которых дифференцируема
функция $F(x)=\int_0^x f(t)dt$,  т.е.  существует предел
$$
\lim_{h\to 0}\frac{F(x+h)-F(x)}{h}=\lim_{h\to 0}\frac1h\int_0^hf(x+t)dt=d_f(x)\quad
(d\text{-точки}). \eqno (3)
$$
Из сходимости во всех $d$-точках следует сходимость в $l-$точках  ($d_f(x)=l_f(x)$), а из
сходимости во всех $l$-точках следует сходимость по норме на всём пространстве
$L_1(\mathbb T)$ (и $C(\mathbb T))$.

Г. Харди доказал ([10],  теорема 253),  что $(C, \alpha)-$ средние рядов Фурье при $\alpha >1$
сходятся во всех $d$-точках и отметил, что $(C, 1)$-средние $\sigma_n(f)$ (см. (2)) могут и
расходиться в $d$-точках, делая ссылку на [6]. Н.К. Бари ([8], гл. I, $\S$49) называет статью
Лебега [11], в которой есть этот результат о $\sigma_n(f)$. Автор не нашёл доказательства
этого факта и получил его из следующего критерия сходимости $\sigma_n(f)$ в $d$-точке
(см. ниже следствие \ref{co4}).

\begin{theorem}\label{th1}
Если $x$ -- $d$-точка функции $f\in L_1(\mathbb T)$, а $$F_x(t)=\frac 1t\int_0^tf(x+u)du,$$ то
$$
\lim_{n\to \infty}(\sigma_n(f;x)+s_n(F_x;0))=d_f(x)+F_x(0)=2d_f(x).
$$
Более того,
$$
\sigma_n(f;x)+s_n(F_x;0)-2d_f(x)=O\(\omega(F_x;\frac {\ln n}{n})+\frac 1n\),
$$
где $\omega$ -- модуль непрерывности $F_x$ на $[-\pi, \pi]$.
\end{theorem}

Кроме того, в настоящей статье доказано общее достаточное условие сходимости средних (1)
во всех $d$-точках в случае $\lambda_{k, n}=\phi(\frac {k}{n+1})$ или, что то же
самое, $\lambda_{k, \epsilon}=\phi (\epsilon k)$ $(\epsilon \searrow 0)$. Для его формулировки
напомним определение винеровской банаховой алгебры:
$$
A(\mathbb R)=\{f:f(x)=\int_{-\infty}^{+\infty}g(t)e^{-itx}dt,\quad
\|f\|_A=\|g\|_{L_1(\mathbb R)}<\infty\}. \eqno(4)
$$

\begin{theorem}\label{th2}
Пусть вариация функции $\phi$ конечна ($\phi\in V$) в некоторой окрестности
нуля, $\phi(0)=1$,  $x\phi(x)\in L_1(\mathbb R)$, $\phi\in A(\mathbb R)$, а $ x\phi'{x}\in
A(\mathbb R)\cap L_1(\mathbb R)$. Тогда во всех $d$-точках любой функции $f\in
L_1(\mathbb T)$
$$
\lim_{n\to
\infty}\sum_{k=-\infty}^{+\infty}\phi\(\frac{k}{n+1}\)\widehat{f}_ke^{ikx}=d_f(x).
$$
\end{theorem}

Более слабое утверждение доказано в [12]. Отметим ещё, что недавно появилась обзорная
статья [13] о винеровских алгебрах.

В \S 2 приведены 4 леммы; доказательство теорем \ref{th1},  \ref{th2} и следствий из них см. в \S 3. В \S 4
приведены новые примеры (методы суммирования Валле-Пуссена, Фейера-Джексона, Рисса, типа
Абеля-Пуассона,  типа Рогозинского-Бернштейна)  и три замечания.

\medskip

\section{Вспомогательные предложения}

Доказательство теоремы \ref{th1}, приведенной во введении, основано на следующем равенстве.

\begin{lemma}\label{le1}
$$
2\sin\frac t2 \Phi'_n(t)=2\Phi_n(t)\cos\frac t2 -D_n(t)\cos\frac t2-\frac 12\cos(n+\frac
12)t+O(\frac 1n),
$$
где
$$
D_n(t)=\frac 12 + \sum_{k=1}^n \cos kt=\frac {\sin (n+\frac 12)t}{2\sin\frac t2} ,
$$
а
$$
\Phi_n(t)=\frac 12+\sum_{k=1}^n (1-\frac {k}{n+1})\cos kt=\frac {\sin^2 (n+1)\frac
t2}{2(n+1)\sin^2\frac t2}
$$
\end{lemma}

{\bf Доказательство.}  Получаем последовательно
\begin{equation*}
  \begin{split}
     &2\sin\frac t2 \Phi'_n(t)\\
     &= \sum_{k=1}^{n+1}k(1-\frac{k}{n+1})\cos(k-\frac12)t -
\sum_{k=1}^{n+1}(k-1)(1-\frac{k-1}{n+1})\cos(k-\frac12)t\\
&=\sum_{k=1}^{n+1}(1-\frac{2k}{n+1}+\frac {1}{n+1})
\cos(k-\frac12)t\\
&=\sum_{k=0}^{n}(1-\frac{2k}{n+1}-\frac {1}{n+1}) \cos(k+\frac12)t,
   \end{split}
\end{equation*}
а значит, это равно и среднему арифметическому последних двух сумм, т.е.,
$$
\cos \frac t2 \sum_{k=1}^{n}(1-\frac {2k}{n+1}) \cos kt +\frac
{1}{2(n+1)}\sum_{k=1}^{n}(\cos(k-\frac 12)t-\cos(k+\frac 12)t)
$$
$$
+ \frac 12(1-\frac
{1}{n+1})\cos \frac t2 -\frac 12(1-\frac {1}{n+1})\cos (n+\frac 12)t$$
$$=2\cos \frac t2 \sum_{k=1}^{n}(1-\frac
{k}{n+1}) \cos kt -\cos \frac t2 \sum_{k=1}^{n} \cos kt
$$
$$
+ \frac 12 \cos \frac t2-\frac 12 \cos (n+ \frac 12)t + O(\frac 1n)
$$
$$
=2\cos \frac t2(\Phi_n(t)-\frac 12)-\cos \frac
t2 (D_n(t)-\frac 12)
$$
$$
+ \frac 12 \cos \frac t2-\frac 12 \cos (n+ \frac 12)t + O(\frac 1n)
$$
$$
2\cos \frac t2 \Phi_n(t)-\cos \frac t2 D_n(t) -\frac 12 \cos (n+ \frac 12)t + O(\frac
1n).
$$
Лемма доказана.

\medskip

\begin{lemma}\label{le2}
Если $\lambda_{0, n}=1 , \lim_{n\to \infty}\lambda_{k, n}=1$ $(k\in \Z)$, $\lim_{|k|\to
\infty}\lambda_{k, n}=0$ $(n\in \N)$ и
$$
\sup_{n}\sum_{k\in \Z}|\Delta\lambda_{k, n}|=\sup_{n}\sum_{k\in \Z}|\lambda_{k, n}-\lambda_
{k+1, n}|<\infty,
$$
ядра $K_n$ абсолютно непрерывны,  а при некотором $\delta\in (0, \pi]$ (в обозначениях
(1))
$$
 \sup_{n}\int_{-\delta}^{\delta} |tK'_n (t)|dt<\infty,
$$
то во всех $d$-точках функции $f\in L_1(\mathbb T)$
$$
 \lim_{n\to \infty}\frac {1}{2\pi}\int_{-\pi}^{\pi}f(x-t)K_n(t)dt=d_f(x).
$$
\end{lemma}
В частном случае чётного ядра $K_n$ и $\delta=\pi$ зто утверждение содержится в
[10] (теорема 71) и [6]. Но доказательство, по сути, не меняется.

Переведём теперь интегральное условие леммы \ref{le2} в условия на коэффициенты
$\{\lambda_{k, n}\}$ ядра $K_n$.

 \begin{lemma}\label{le3}

Если $K(t)=\sum_{k\in Z}\lambda_k e^{ikt}$ и $\lim_{|k|\to \infty}k\lambda _k=0$, то
$$
\int_{-\pi}^{\pi}|tK'(t)|dt\le \frac \pi 2\int_{-\pi}^{\pi}|K(t)|dt+\frac \pi
2\int_{-\pi}^\pi |\sum_{k\in \mathbb Z}k\Delta\lambda_k e^{ikt}|dt,
$$
при условии, что правая часть конечна.
\end{lemma}

{\bf Доказательство.}

В силу неравенства $\frac 2\pi |u|\le |\sin u|$ при $|u|\le \frac \pi 2$
$$
\int_{-\pi}^{\pi}|tK'(t)|dt\le\frac \pi 2\int_{-\pi}^{\pi}|2\sin \frac t2 K'(t)|dt=
\frac \pi 2 \int_{-\pi}^{\pi}|2\sin \frac t2 \sum_{k\in \mathbb Z}k\lambda_k e^{ikt}|dt=
$$
$$
\frac \pi 2\int_{-\pi}^{\pi}|\sum_{k\in \mathbb Z}k\lambda_k( e^{ikt} - e^{i(k+1)t}|dt=
\frac \pi 2\int_{-\pi}^{\pi}|\sum_{k\in \mathbb Z}\Delta (k\lambda_k)e^{ikt}|dt.
$$
Осталось учесть, что
$$\Delta(k\lambda_k)=k\lambda_k-(k+1)\lambda_{k+1}=k\Delta\lambda_k-\lambda_{k+1}.$$

\begin{lemma}\label{le4}
Если $f\in A(\mathbb R)\cap L_1(\mathbb R), $ то преобразование Фурье $\widehat {f}\in
L_1$.
\end{lemma}

{\bf Доказательство.}

Если $g\in L_1$, а преобразование Фурье $$\widehat {g}(y)=\frac {1}{\sqrt
2\pi}\int_{-\infty}^{\infty}g(x)e^{-ixy}dx, $$ то $ \check {g}(y)=\widehat
{g}(-y)$ -- обратное преобразование Фурье. Если и $\widehat g \in L_1$, то, как известно,
почти всюду
$$
g=\check {\widehat g}=\widehat {\check g}.
$$
Докажем, что если
$$
F=\widehat {F_1}, F_1\in L_1, F_1=\check {F_2}, F_2\in L_1,
$$
то $F=F_2$ почти всюду.

В силу формулы умножения, применённой дважды (считаем $g$ и $\widehat g\in L_1$),
$$
\int _\mathbb R Fg=\int_\mathbb R \hat {F_1}g=\int _\mathbb R F_1\hat g=\int _\mathbb R
\check {F_2}\widehat g=\int _\mathbb R F_2g.
$$
Так что для всех таких функций $g$
$$ \quad \int _\mathbb R (F-F_2)g=0,$$
где $F\in C(\mathbb R)$, а $F_2\in L_1. $

Докажем, что $F=F_2$ почти всюду на любом отрезке $[a, b]$. Для этого достаточно
доказать, что $\int_a^b(F-F_2)=0$, так как тогда из равенства $\int_a^x(F-F_2)=0$  сразу
следует, что производная $F-F_2=0 $ почти всюду.

Для аппроксимации функции $g_0$, которая равна единице на $[a, b]$ и нулю вне
$[a, b]$, возьмём непрерывную функцию $g_n$, которая равна единице на $[a+\frac 1n, b-\frac
1n], g_n(a)=g_n(b)=0$ и линейная на $[a, a+\frac 1n]$ и $[b-\frac 1n, b]. $ Очевидно, что
$g_n$ и  $\widehat g_n\in L_1$, поэтому $\int _a^b(F-F_2)g_n=0. $

Кроме того,
$$
\bigg|\int_a^b F\cdot (g-g_n)\bigg|\le \|F\|_{C_{[a, b]}}\int_a^b|g-g_n|\le \frac 2n
\|F\|_{C_{[a, b]}}
$$
и
$$
\bigg|\int_a^b F_2\cdot (g-g_n)\bigg|\le \int_a^b |F_2||(g-g_n)|
$$
 и учитывая ещё, что $|g-g_n|\le 1$, можно перейти к пределу при $n\to \infty$ в силу
 теоремы Лебега о мажорируемой сходимости. Получаем, что $\int _a^b (F-F_2)=0,$ и
 $F=F_2$ почти всюду.

 При условиях леммы \ref{le4} берём $F=\widehat {f}, f\in L_1. $ Тогда $F_1=f\in A(\mathbb R). $ Поэтому
 $F_1=\check F_2, F_2\in L_1. $ По доказанному $\widehat {f}=F=F_2\in L_1. $

 Лемма \ref{le4} доказана.

\medskip

\section{Доказательство теорем. Следствия}

\medskip

{\bf Доказательство теоремы \ref{th1}, приведенной
во введении.}

Пусть сначала $f_0\in L_1(\mathbb T)$, $d_f(0)=0$ и $$F_0(t)=\frac 1t\int_0^tf(u)du, \quad
F_0(0)=0,  \quad F_0\in C[-\pi, \pi]).$$ Тогда, применяя интегрирование по частям, получаем
$$
\sigma_n(f_0;0)=\frac 1\pi \int _{-\pi}^{\pi} f_0(t)\Phi_n(t)dt= \frac 1\pi \Big[\frac
{\sin^2 (n+1)\frac t2}{2(n+1)\sin^2 \frac t2}\cdot \int_0^t f_0(u)du\Big]_{-\pi}^{\pi}+
$$
$$
 \frac 1\pi\int_{-\pi}^{\pi}F_1(t)2\sin \frac t2 \Phi'_n(t)dt,
$$
где $F_1(t)=F_0(t)\frac {t}{2\sin \frac t2}$. Внеинтегральный член есть $O(\frac 1n)$.

Применяем лемму \ref{le1}:
$$
\sigma_n(f_0;0)=\frac 2\pi\int_{-\pi}^{\pi}F_1(t)\cos \frac t2\Phi_n(t)dt-\frac
1\pi\int_{-\pi}^{\pi}F_1(t)\cos \frac t2 D_n(t)dt$$
$$
-\frac
12\int_{-\pi}^{\pi}F_1(t)\cos(n+\frac12)tdt+O(\frac 1n).
$$
Теперь от $F_1$ вернёмся к $F_0$ с оценкой погрешности.

Если $h\in {\rm Lip}1$ на $[-\pi, \pi]$, то продолжая $F_0$ и $h$ нулём на $\mathbb R\setminus
[-\pi, \pi] $, имеем при $\lambda \ne 0$
$$
\int_{-\infty}^{\infty}F _0(t)h(t)e^{i\lambda t}dt=O\(\omega(F_0;\frac
{1}{|\lambda|})+\frac {1}{|\lambda|}\)  \eqno (5)
$$
(доказательство приведено ниже).

Очевидно, что функция $h(t)=(\frac 12 tctg \frac t2-1)\frac {1}{2\sin^2\frac t2}$
ограничена на $[-\pi, \pi]$, как и её производная. Поэтому, учитывая ещё
 ограниченность $F_0$ и $\sin$, получаем
 \begin{equation*}
   \begin{split}
 \int_{-\pi}^{\pi}(F_1(t)\cos &\frac t2-F_0(t))\Phi_n(t)dt\\
&=\frac {1}{n+1}\int
 _{-\pi}^{\pi}F_0(t)h(t)\sin^2(n+1)\frac t2 dt=O\(\frac 1n\),
    \end{split}
 \end{equation*}
  \begin{equation*}
   \begin{split}
\int_{-\pi}^{\pi}(F_1(t)\cos &\frac t2-F_0(t))D_n(t)dt\\
&=\int
 _{-\pi}^{\pi}F_0(t)h(t)\sin \frac t2 \sin (n+\frac 12)tdt=O\(\omega(F_0;\frac 1n)
+\frac 1n\)
    \end{split}
 \end{equation*}
и
$$
\int_{-\pi}^{\pi}(F_1(t)-F_0(t))\cos (n+\frac 12)tdt=O\(\omega(F_0;\frac 1n) +\frac 1n\).
$$
Таким образом,
$$
\sigma_n(f_0;0)=2\sigma_n(F_0;0)-s_n(F_0;0)+O(\omega(F_0;\frac 1n) +\frac 1n)
$$
Оценим скорость стремления к нулю $\sigma_n(F_0;0)=\sigma_n(F_0;0)-F_0(0)$ через модуль
непрерывности $F_0$ на $[-\pi, \pi]$.
$$
|\sigma_n(F_0;0)|=\frac 1\pi\bigg|\int_{-\pi}^{\pi}(F_0(t)-F_0(0))\Phi_n(t)dt\bigg|\le
\frac1\pi\int_{-\pi}^{\pi}\omega(F_0;|t|)\Phi_n(t)dt.
$$
Известно, что при $\lambda>0$ $$\omega(f;\lambda|u|)\le (\lambda+1)\omega(f;|u|). $$
Поэтому
$$
|\sigma_n(F_0;0)|\le \omega(F_0;\frac {\ln(n+1)}{n+1})\Big (\frac 1\pi
\int_{-\pi}^{\pi}\frac {(n+1)|t|}{\ln(n+1)}\cdot \frac {\sin^2(n+1)\frac
t2}{2\sin^2\frac t2}\cdot \frac {dt}{n+1} +1\Big).
$$
После применения в знаменателе неравенства $|u|\le \frac \pi 2|\sin u|$ при $|u|\le
\frac \pi 2$ первое слагаемое в скобках не больше
$$
\frac {\pi}{\ln(n+1)}\int_0^\pi\frac {\sin^2(n+1)\frac t2}{t}dt=\frac
{\pi}{\ln(n+1)}\int_{0}^{\frac 12\pi(n+1)}\frac {\sin^2u}{u}du=O(1).
$$
Так что
$$
|F_0(0)-\sigma_n(F_0;0)|\le c\omega\(F_0;\frac {\ln (n+1)}{n+1}\).\eqno (6)
$$
Следовательно,
$$
\sigma_n(f_0;0)+s_n(F_0;0)=O\(\omega(F_0;\frac {\ln (n+1)}{n+1})\)+O\(\frac 1n\).
$$
Пусть теперь $x$ -- произвольная $d$-точка функции $f$.

Применим доказанное соотношение к функции $f_0(t)=f(x+t)-d_f(x)$. При этом
$$
F_0(t)=\frac 1t\int_0^t f(x+u)du -d_f(x).
$$
Тогда
$$
\sigma_n(f_0;x)+s_n(F_x;0)-2d_f(x)=O\(\omega(F_x;\frac {\ln (n+1)}{n+1}\)+\frac 1n).
$$
Осталось доказать неравенство (5).

После простых преобразований ($\|g\|_\infty=ess\sup_{[\pi, \pi]}|g(t)|$)
$$
\bigg|\int_{-\infty}^{\infty}F_0(t)h(t)e^{i\lambda t}dt\bigg|=\frac
12\bigg|\int_{-\infty}^{\infty}(F_0(t)h(t)-F_0(t+\frac \pi\lambda))h(t+\frac \pi\lambda))dt\bigg|=
$$
$$
\frac 12\bigg|\int_{-\infty}^{\infty}(F_0(t)-F_0(t+\frac {\pi}{\lambda})h(t)e^{i\lambda
t}dt+\int_{-\infty}^{\infty}F_0(t+\frac {\pi} {\lambda})(h(t)-h(t+\frac
{\pi}{\lambda}))e^{i\lambda t}dt\bigg|
$$
$$
\le\frac 12\|h\|_\infty\int_{-\infty}^{\infty}|F_0(t)-F_0(t+\frac {\pi}{\lambda})|dt$$
$$+\frac12\|F_0\|_\infty\int_{-\infty}^{\infty}|h(t)-h(t+\frac {\pi}{\lambda})|dt.
$$
При $\lambda >0$, например, первый интеграл равен
$$
\int_{-\pi-\frac {\pi}{\lambda}}^\pi|F_0(t+\frac
{\pi}{\lambda})|dt+\int_{-\pi}^{\pi-\frac {\pi}{\lambda}}|F_0(t)-F_0(t+\frac
{\pi}{\lambda})|dt+\int_{\pi-\frac {pi}{\lambda}}^\pi|F_0(t)|dt
$$
$$\le
2\|F_0\|_\infty\cdot\frac {\pi}{\lambda}+2\pi\omega(F_0;\frac {\pi}{\lambda}).
$$
При такой же оценке второго интеграла нужно ещё учесть, что $\omega(g; \frac
{\pi}{\lambda})\le\|g'\|_\infty \frac {\pi}{\lambda}. $

Неравенство (5) доказано, а с ним и теорема \ref{th1}.

Приведём несколько следствий.

\begin{corollary}\label{co1}
 Если $x$ -- $l$-точка функции $f\in L_1(T)$, то ряд Фурье
функции $F_x(t)=\frac 1t\int_0^tf(x+u)du$ в нуле сходится.
\end{corollary}

Для доказательства достаточно применить теорему \ref{th1} и теорему Лебега, упомянутую во
введении.

\begin{corollary}\label{co2}
Для того чтобы $\lim_{n\to \infty}\sigma_n(f;x)=d_f(x)$, необходимо и достаточно, чтобы
сходился в нуле ряд Фурье непрерывной функции $F_x$.
\end{corollary}
\medskip

\begin{corollary}\label{co3}
$(C, \alpha)$-средние рядов Фурье при $\alpha >1$ сходятся во всех $d-$точках.
\end{corollary}

Для доказательства на основании теоремы \ref{th1} достаточно применить к предельному соотношению
теоремы $(C, \epsilon)$-средние при $\epsilon >0$. Получим
$$
\lim (\sigma_n^{1+\epsilon}(f;x)+\sigma_n^\epsilon(F_x;0))=2d_f(x).
$$
Но $F_x\in C[-\pi. \pi]$ и поэтому $\sigma_n^\epsilon (F_x;0)\to F_x(0)=d_f(x). $

Этот же приём применим и в более общей ситуации.

\begin{corollary}\label{co4}
Существует функция $f\in L_1(\mathbb T)$, у которой $x=0$ является $d$-точкой, а
$\sigma_n(f;0)= \sigma^1_n(f;0)$ при $n\to \infty$ расходится.
\end{corollary}

{\bf Доказательство.} Применяем следствие \ref{co2}.
Есть много разных примеров расходящихся рядов Фурье непрерывных функций (см. [5-8]).  Нам
нужна четная функция $F_0$ вида
$$
tF_0(t)=\int_0^tf_0(u)\, du, \quad f_0\in L_1({\mathbb T}),
$$
с расходящимся в нуле рядом Фурье.

Воспользуемся примером из [6], \S4.12 (см. также [8], гл. I, \S 46).

При $n_k=3^{k^3}\ (k\geq0)$ и $a_k=\frac1{k^{3/2}}\ (k\in\mathbb N)$
$$
F_0(t)=a_k\sin n_k|t|\quad \Big(\frac\pi{n_k}\leq|t|\leq\frac\pi{n_{k-1}}\Big).
$$
Поскольку $\frac{n_k}{n_{k-1}}$ -- нечетное число,  а
$\lim_{k\to\infty}a_k\ln\frac{n_k}{n_{k-1}}=\infty$,  то
$$
\lim_{k\to\infty}s_{n_k}(F_0;0)=\infty.
$$
Функция $F_0$ не только непрерывна $(F_0(0)=F_0(\pm\pi)=0)$,  но и имеет конечные
односторонние производные в точках $\pm\frac\pi{n_k}\ (k\geq0)$.  При $t>0$
$$
tF_0(t)=-\int_t^\pi (F_0(u)+uF'_0(u))\, du.
$$
Но
$$
\int_0^\pi|uF'_0(u)|\, du=\sum_{k=1}^\infty\int_{\frac\pi{n_k}}^{\frac\pi{n_{k-1}}}
ua_kn_k|\cos n_ku|\, du$$
$$
\leq\sum_{k=1}^\infty
a_kn_k\frac12\Big(\frac\pi{n_{k-1}}\Big)^2<\infty.
$$
Следствие \ref{co4}  доказано.

Переходим к {\bf доказательству теоремы \ref{th2}}.

Применяем лемму~\ref{le2} при $\delta=\pi$ и лемму \ref{le3}.  Функция $\phi\in
C(\mathbb R)\cap V(\mathbb R)$, так как $A(\mathbb R)\subset C(\mathbb R)$ и $\phi'\in
L_1(\mathbb R)$.

Проверим, что $\phi(x)=o(\frac 1x)$  при $|x|\to \infty$. Действительно,
$$
\phi(x)=-\int_x^{+\infty\cdot {\rm sign} x}\phi'(t)dt=-\int_x^{+\infty\cdot {\rm sign} x}t\cdot
\frac 1t\phi'(t)dt,
$$
а из условия $x\phi'(x)\in L_1(R)$ следует, что
$$
|\phi(x)|\le \frac {1}{|x|}\int_{|x|}^\infty t(|\phi'(t)+|\phi'(-t)|)dt=o\(\frac 1x\).
\eqno(7)
$$
Для оценки сверху интеграла от модуля ряда Фурье есть много
результатов (см., например, [7], п. 7.2).  Воспользуемся следующей теоремой из [12] (теорема
8): если $\sum \psi(k)e_k$ -- ряд Фурье функции $\Psi$, то
$$
\frac 1{2\pi}\int_{-\pi}^\pi|\Psi(t)|dt=\min_{\psi_c}\|\psi_c\|_A,
$$
где $\psi_c$ -- любая непрерывная функция с условием $\psi_c(k)=\psi(k)$ $(k\in \Z)$.

Учитывая, что при $\lambda\ne 0$  выполняется $\|f(\lambda\cdot)\|_A=\|f(\cdot)\|_A$, получаем
$$
\frac {1}{2\pi}\int_{-\pi}^\pi|K_n(t)|dt=\frac {1}{2\pi}\int_{-\pi}^\pi\bigg|\sum_k
\phi\bigg(\frac k{n+1}\bigg)e^{ikt}\bigg|dt
$$
$$
\le
\|\phi((n+1)^{-1}{\cdot})\|_A=\|\phi(\cdot)\|_A=\|\phi\|_A.
$$
Для такой же оценки второго интеграла в лемме \ref{le3} положим
$$
k\Delta\lambda_{k, n}=\phi_n(\frac k{n+1}), \quad \phi_n(x)=nx(\phi(x)-\phi(x+\frac1n)).
$$
Нужно доказать, что
$$
\phi_n(x)=\int_{-\infty}^\infty g_n(y)e^{-ixy}dy, \quad
\sup_n||\phi_n||_A=\sup||g_n||_{L_1}<\infty. \eqno (8)
$$
Из того, что $\phi\in A(\mathbb R)$, следует, что
$$
\phi(x)=\int_{-\infty}^{\infty}g(y)e^{-ixy}dy,\quad g\in L_1(\mathbb R).
$$
Но тогда
$$
\phi_n(x)=nx\int_{-\infty}^{\infty}g(y)(1-e^{-\frac{iy}{n+1}})e^{-ixy}dy. \eqno (9)
$$
По формуле обращения ($g$ и ее преобразование Фурье принадлежат $L_1$)
$$
g(y)=\frac{1}{2\pi}\int_{-\infty}^{\infty}\phi(x)e^{ixy}dx
$$
почти всюду, а можно считать и всюду на $\mathbb R$. В силу леммы Римана-Лебега
$\lim_{|y|\to \infty}g(y)=0$.  Учитывая теперь, что по условию теоремы $x\phi(x)\in
L_1(\mathbb R)$,  получаем
$$
g'(y)=\frac{i}{2\pi}\int_{-\infty}^{\infty}x\phi(x)e^{ixy}dx,
$$
а после интегрирования по частям (см. ещё (7))
$$
yg'(y)=-\frac {1}{2\pi}\int_{-\infty}^{\infty}(\phi(x)+x\phi'(x))e^{ixy}dx. \eqno (10)
$$
После интегрирования по частям в (9) имеем
$$
\phi_n(x)=-in\int_{-\infty}^{\infty}\bigg[g'(y)(1-e^{-\frac {iy}{n+1}})+g(y)\frac in
e^{-\frac {iy}{n+1}}\bigg]e^{ixy}dy.
$$
Учитывая , что $|e^{ix}-1|\le |x|$, приходим к выводу (см. ещё (8)):
$$
\|\phi_n\|_A\le\int_{-\infty}^{\infty}(|yg'(y)|+|g(y)|)dy.
$$
Осталось доказать, что  $yg'(y)\in L_1(\mathbb R). $

Воспользуемся леммой \ref{le4}. Достаточно применить эту лемму к функции $\phi(x)+x\phi'(x)$
(см. (10)).  Теорема \ref{th2} доказана.

Применяя признак Зигмунда ([5], т. I, гл. VI, (3.6), cм. также [7], 6.4.3), получаем

\begin{corollary}\label{co5}
Если $\phi\in C(\mathbb R)\cap C^1(\mathbb R\setminus \{0\})$, $\phi(0)=1$, ${\rm supp\,} \phi\subset
[-1, 1], $ а $\phi(x)$  и $x\phi'(x)\in V\cap {\rm Lip}\, \epsilon, \epsilon >0, $ то
$$
\lim_{n\to \infty} \sum _{|k|\le n}\phi\(\frac {k}{n+1}\)\widehat {f}_ke^{ikx}=d_f(x).
$$
\end{corollary}

\section{Примеры. Замечания}

\medskip

{\bf Пример 1.}

Ядро Фейера-Джексона
$K_n(t)=\gamma_{s, n}D_n^s(t)$, $\int_{-\pi}^{\pi}K_n=2\pi,$ удовлетворяет условиям леммы \ref{le2}
при любом натуральном $s\ge 3$.

Эта лемма применима и в случае, когда при некотором $\delta\in (0, \pi]$ производная ядра
$K_n$ при любом $n\in N$ сохраняет знаки на $[0, \delta]$ и $[-\delta, 0]$. Тогда на
$[0, \delta]$, например,
$$
\int_0^\delta t|K'_n(t)|dt=\bigg|\int_0^\delta tK_n(t)dt\bigg|=\bigg|\delta K_n(\delta) -\int_0^\delta
K_n(t)dt\bigg|
$$
$$
\le c\sup_n \sum_{k\in \Z}|\Delta\lambda_{k, n}|+\sup_n\int_0^\delta|K_n(t)|dt.
$$
А если к тому же есть сходимость на $L_1(\mathbb T)$, то нормы операторов $\Lambda_n$,
равные $ \frac {1}{2\pi}\int_{-\pi}^{\pi}|K_n(t)|dt$, ограничены.

\medskip

{\bf Пример 2.}

Средние Валле-Пуссена (см., например, [10], 4.17)
$$
\Lambda_n(f;x)=\gamma_n\int_{-\pi}^{\pi}f(x-t)\cos^n\frac t2 lt,
\quad\Lambda_n(1:x)\equiv 1
$$
сходятся во всех $d-$точках.

\medskip

{\bf Пример 3.}

Функция $\phi(x)=e^{-|x|^\alpha}$ удовлетворяет условиям теоремы \ref{th2} при любом $\alpha
>0$. Ранее сходимость в $d-$точках была известна лишь при $\alpha=1$ (метод Абеля-Пуассона) и
$\alpha=2$ (см. [10], п.3 приложения II).

\medskip

{\bf Пример 4.}

Рассмотрим $\varphi(x)=(1-|x|^\alpha)_+^\beta, \ \alpha>0$  и $\beta>0$ (средние Рисса).
Сходимость в $d$-точках имеет место при любом $\alpha>0$,  но только при $\beta>1$.

При $\beta>1$ применяем следствие \ref{co5}. При $\beta=1$ то же следствие применяем к разности
$$
(1-|x|^\alpha)_++\alpha(1-|x|)_+.
$$
Но для $(C, 1)$ нет сходимости во всех $d$-точках для всех $f\in L_1({\mathbb T})$
(см. следствие \ref{co4}).  Следовательно,  и здесь нет.  А если бы такая сходимость была при
$\beta<1$,  то она была бы и при $\beta=1$ (см.  соответствующую теорему для произвольных
числовых рядов в [10], теорема~58).

\medskip

{\bf Пример 5} (метод типа Рогозинского-Бернштейна).

Рогозинский изучал равномерную сходимость при $n\to\infty$ средних
$$
\frac12\Big(s_n \Big(f;x+\frac\pi{2n}\Big)+s_n\Big( f;x-\frac\pi{2 n}\Big)\Big)\quad
\bigg(s_n (f)=\sum_{k=-n}^n\hat f_ke_k\bigg),
$$
а С.Н.  Бернштейн -- близких по идее средних
$$
\frac12\Big(s_n(f;x)+s_n\Big(f;x+\frac\pi n\Big)\Big).
$$
При сравнении скорости сходимости (по норме) их аппроксимативные свойства оказались
разными: в первом случае точный порядок приближения $\omega_2 (f;\frac1n)$ (модуль
гладкости второго порядка), во втором -- $\omega(f;\frac 1n)$ (см., например, [7], 8.5.1).

Рассмотрим общие средние типа Рогозинского--Бернштейна
$$
\int_{-\infty}^\infty s_n(f;x-\varepsilon\gamma
t)\, d\mu(t)=\sum_{|k|\leq\frac1\varepsilon}\varphi(\varepsilon k)\hat f_ke^{ikx},
$$
где при $|x|\leq1$
$$
\varphi(x)=\int_{-\infty}^\infty e^{-i\gamma xt}d\mu(t), \quad \gamma\in\mathbb R, \
\varphi(x)=0 \ (|x|\geq1), \ \varphi(0)=1,
$$
а $\mu$~ -- конечная на $\mathbb R$ комплекснозначная мера (см. пример~III после 8.1.4 в
[7]).

В силу теоремы \ref{th2},  если еще $\varphi'(\pm1)=0$ и
$$
\int_{-\infty}^\infty|t|^2|d\mu(t)|<\infty,
$$
то имеем сходимость во всех $d$-точках.  Если же не выполняется условие
$\varphi'(\pm1)=0$,  то применяем теорему к разности
$$
\varphi(x)-\frac12(\varphi'(1)-\varphi'(-1))(1-|x|)_+-\frac12(\varphi'(1)
+\varphi'(-1))(1-|x|)_+{\rm sign}\, x.
$$
Для четных функций $f$ нужна информация о поведении средних арифметических $\sigma_n(f)$
(см. следствие \ref{co4}),  а для нечетных функций известно,  что если $0$~ -- точка Лебега,  то для
всех $f\in L_1({\mathbb T})$ с условием $\int_{\to 0}^\pi\frac{f(t)}t\, dt=\infty$, $
\lim_{n\to\infty}\,\sigma_n(f, 0)=\infty$ (см.  [5], т. ~I,  гл. ~III,  (3.20)).  А средние
Рогозинского (финитная функция $\varphi(x)=\cos\frac{\pi x}2, \ |x|\leq1$) и средние
Бернштейна (финитная функция $\varphi(x)=\cos^2\frac{\pi x}2+\frac i2\sin\pi x, \ |x|\leq
1$) сходятся во всех точках Лебега,  так как $\sup_{|y|\geq x}|\hat\varphi(y)|\in
L_1(\mathbb R_+)$ (см. 8.1.3 в [7]),  но не всегда~ -- в $d$-точках.  Если же применить к
ним еще средние Рисса при $\alpha>0,  \beta>0$,  то получим сходимость во всех $d$-точках.

\medskip

{\bf Замечание 1.}

 В дополнение к лемме \ref{le1} приведём ещё несколько равенств, которые можно применить, например, к
 сопряжённым рядам Фурье.

 При
 $$
 \tilde {D}_n(t)=\sum_{k=1}^n \sin kt,\quad  \tilde {\Phi}_n(t)=\frac
 {1}{n+1}\sum_{k=1}^n \tilde {D}_k(t)
 $$
  1)  \begin{equation*}
   \begin{split}
 &2\sin\frac t2 \tilde {\Phi}'_n(t)\\
&=-2\cos \frac t2 \tilde {\Phi}(t)+\cos\frac t2\tilde
{D}_n(t)-\frac 12\sin \frac t2 +\frac 12\sin (n+\frac 12)t+O(\frac 1n).
    \end{split}
 \end{equation*}
  2)  \begin{equation*}
   \begin{split}
(1-e^{-it})&\sum_{k=1}^n(1-\frac {k}{n+1})ke^{ikt}\\
&=(1+\frac
{1}{n+1})\sum_{k=0}^{n}e^{ikt}-2\sum_{k=0}^{n}(1-\frac {k}{n+1})e^{ikt},
    \end{split}
 \end{equation*}
  3)  \begin{equation*}
   \begin{split}
(1-e^{it})&\sum_{k=1}^n(1-\frac {k}{n+1})ke^{ikt}\\
&=-(1-\frac
{1}{n+1})\sum_{k=1}^{n+1}e^{ikt}+2\sum_{k=1}^{n}(1-\frac {k}{n+1})e^{ikt}.
    \end{split}
 \end{equation*}

\medskip

{\bf Замечание 2} (о признаке Салема равномерной сходимости рядов Фурье).

Применим теорему \ref{th1} к периодическим функциям из $C(\mathbb T)$ $(d_f(x)=f(x))$.  Так как
$$
F_x(t)=\frac 1t\int _0^tf(x+u)du=\int_0^1f(x+tu)du,
$$
то
\begin{equation*}
   \begin{split}
\omega(F_x;h)&=\sup_{-\pi\le t\le t+\delta\le t+h\le
\pi}|F_x(t)-F_x(t+\delta)|\\
&\le\omega(f;h)=\sup _{0<\delta\le h, t\in T}|f(t)-f(t+\delta)|.
    \end{split}
 \end{equation*}
Получаем, учитывая ещё (6),
$$
\sup_{x\in \mathbb T}|f(x)-s_n(F_x;0)|=O\(\omega(f;\frac {\ln n}{n})+\frac 1n\).
$$
Это некоторая связь между скоростью сходимости ряда Фурье функции $f$ и
проинтегрированной функции.

А в силу признака Салема (см. [5], гл. IV, \S7), если $f\in C(\mathbb T)$ и $\hat {f}_0=0, $ а
$F(x)=\int _0^x f$, то при $n\to \infty$
$$
\|F-s_n(F)\|_\infty =o(\frac 1n) \Longrightarrow  \|f-s_n(f)\|_\infty=o(1).
$$
На самом деле, верно и обратное утверждение, т.е., имеем критерий равномерной сходимости
рядов Фурье (сформулирован в [7], 2.5.12). Для доказательства этого достаточно применить
известное неравенство Бора-Бернштейна (см., например., 5.5.2 в [7] при $r=1$):
$$
\|F-s_n(F)\|_\infty\le \frac {\pi}{2n+2}\|f-s_n(f)\|_\infty.
$$
Если же придерживаться доказательства Салема (см. там же в [8]), то получаем следующий
критерий сходимости в точке ряда Фурье непрерывной функции:
$$
f(x)-s_n(f;x)=o(1) \iff F(x+\frac{\pi}{2n})-s_n(F;x+\frac {\pi}{2n})=o(\frac 1n).
$$
Действительно, равномерно по $x\in T$ (см. [8], гл. IV, (7. 4)) при $h=\frac {\pi}{2n}$
$$
\frac 12(s_n(F;x+\frac {\pi}{2n})-s_n(F;x)-\frac {\pi}{2n}f(x)=\frac 1n (s_n(f;x)-f(x))+
o(\frac 1n). \eqno (11)
$$
Учтём теперь, что при $f(x)=F'(x)\in C(\mathbb T)$
$$
\frac {\pi}{2n}f(x)=\frac 12(F(x+\frac {\pi}{2n})-F(x-\frac {\pi}{2n}))+o(\frac 1n)\eqno
(12)
$$
и
$$
\omega_2(F;h)\le h\omega (F';h)=h\omega (f;h) =o(h).
$$
Но тогда
$$
F(x)-\frac 12(F(x+\frac {\pi}{2n})+F(x-\frac {\pi}{2n}))=o(\frac 1n)\eqno (13)
$$
и (см. ещё точный порядок приближения в примере 5)
$$
\frac 12 (s_n(F;x+\frac {\pi}{2n})+s_n(F;x-\frac {\pi}{2n}))-F(x)=o(\frac 1n).\eqno (14)
$$
Складывая равенства (11)-(14) и умножая сумму на (-1), получаем
$$
F(x+\frac {\pi}{2n})-s_n(F;x+\frac {\pi}{2n})=\frac 1n (f(x)-s_n(f;x))+o(\frac 1n)
$$
(здесь можно $+\frac {\pi}{2n} $ заменить на -$\frac {\pi}{2n}$.)

\medskip

{\bf Замечание 3} (о необходимом условии суммируемости).

Для сходимости на всём пространстве $L_1(\mathbb T)$ средних
$$
\Lambda_n(x)=\sum_{|k|\le n}\lambda_{k, n}\hat {f}_k e^{ikx}
$$
(тем более, для сходимости в $l$-точках и $d$-точках) необходима ограниченность норм
операторов $\Lambda_n$:
$$
\sup_n \int_{-\pi}^\pi|\sum_{-n}^n \lambda_{k, n}e^{ikt}|dt<\infty.
$$
Сидон доказал, что при некотором числе $c>0$ при всех $n$ и $\lambda_{k, n}$
$$
 \int_{-\pi}^\pi|\sum_{-n}^n \lambda_{k, n}e^{ikt}|dt\ge c\sum_{-n}^n\frac
 {|\lambda_{k, n}|}{n-|k|+1}.
$$
Приведём это доказательство, основанное на известном неравенстве Харди-Литтльвуда. Имеем
\begin{equation*}
  \begin{split}
     \int_{-\pi}^\pi|\sum_{-n}^n \lambda_{k}e^{ikt}|dt&=\int_{-\pi}^\pi|\sum_{0}^{2n}
 \lambda_{k-n}e^{ikt}|dt\\
 &\ge c_1\sum_{k=0}^n\frac
 {|\lambda_{k-n}|}{k+1}=c_1\sum_{k=-n}^{0}\frac  {|\lambda_{k, n}|}{n-|k|+1}.
   \end{split}
\end{equation*}
Отсюда следует, что
\begin{equation*}
  \begin{split}
     \int_{-\pi}^\pi|\sum_{-n}^n \lambda_{k}e^{ikt}|dt&=\int_{-\pi}^\pi|\sum_{-n}^n
\lambda_{-k}e^{ikt}|dt\\
&\ge c_1\sum_{k=-n}^{0}\frac  {|\lambda_{-k, n}|}{n-|k|+1}=
\sum_{k=0}^{n}\frac  {|\lambda_{k, n}|}{n-|k|+1}.
   \end{split}
\end{equation*}
Соединяя эти два неравенства, получаем
$$
\int_{-\pi}^\pi|\sum_{-n}^n \lambda_{k, n}e^{ikt}|dt\ge \frac 12 c_1\sum_{-n}^n\frac
 {|\lambda_{k, n}|}{n-|k|+1},
$$
где $c_1$ -- константа из неравенства Харди-Литтльвуда.

Но уже получено  более сильное неравенство ([14]):
$$
\int_{-\pi}^\pi|\sum_{1}^{\infty} \lambda_{k}e^{ikt}|dt\ge
c_2\sum_{s=1}^{\infty}\left(\sum_{2^{s-1}\le \nu
<2^s}\frac{|\lambda_{\nu}|^2}{\nu}\right)^{\frac 12}.
$$
Повторяя предыдущие рассуждения, можно получить более сильное необходимое условие для
сходимости по норме в $L_1(\mathbb T)$ и $C(\mathbb T)$.

Заметим теперь, что интегральное условие леммы \ref{le2} не является необходимым для сходимости в
$d$-точках. Для доказательства этого рассмотрим метод суммирования Лебега
$$
\Lambda_n(f;x)=\sum_{-\infty}^\infty\frac {\sin k\epsilon_n}{k\epsilon_n}\hat {f}_k e^{ikx}
\quad (\epsilon_n \to 0).
$$
(см. [8]). Ядро
$$
K_n(t)=\sum_{-\infty}^{\infty}\frac {\sin k\epsilon_n}{k\epsilon_n}e^{ikt}
$$
неотрицательно, так как (воспользуемся ещё раз теоремой 8 из [12])
\begin{equation*}
\begin{split}
1=\frac{1}{2\pi}\int_{-\pi}^{\pi}K_n(t)dt&\le
\frac{1}{2\pi}\int_{-\pi}^{\pi}|K_n(t)|dt\\
&\le\bigg\Vert\frac {\sin x\epsilon_n}{x\epsilon_n}\bigg\Vert_A=\bigg\Vert\frac {\sin
x}{x}\bigg\Vert_A=\bigg\Vert\frac 12 \int_{-1}^1 e^{-ixy}dy\bigg\Vert_A=1.
   \end{split}
\end{equation*}

Суммируемость имеем не только на $L_1(\mathbb T)$, но и в $d$-точках, так как
$$
\lim_{n\to \infty}\Lambda_n(f;x)=\lim_{n\to\infty}\frac
{F(x+\epsilon_n)-F(x-\epsilon_n)}{2\epsilon_n}=d_f(x).
$$
Продифференцированный ряд $K_n$ не является рядом Фурье, так как его коэффициенты не
стремятся к нулю. Поэтому ядро не может быть абсолютно непрерывным. Кроме того,  этот ряд в
силу теоремы Кантора-Лебега (см. [5], [8]) расходится почти всюду.

Укажем теперь простое необходимое условие.

Если  средние рядов Фурье, определяемые ядром $K_n$,  сходятся в $d$-точке
$$
\lim_{h\to 0} \frac 1h\int_0^h g(t)dt=0,\quad \lim\sup_{h\to 0}\frac 1h\int_0^h|g(t)|dt>0
$$
($d_g(0)=0$ и $0$ не является $l$-точкой), то при любом $\delta\in (0, \pi]$
$$
\lim_{n\to \infty}\int_0^{\delta}g(t)K_n(t)dt=0.
$$
В качестве такой функции $g$ возьмём, например, следующую. Пусть
$\{x_s\}_0^\infty$-положительная последовательность, убывающая к нулю, а $x_0=\pi$.
Полагаем $g(t)=(-1)^s$ при $t\in (x_{s+1}, x_s]$ и $g(t)=0$ при $t\in [-\pi, 0]$.Легко
проверить, что $d_g(0)=0$ в том и только в том случае, когда $\lim \frac
{x_{s+1}}{x_s}=1$. Следовательно, для всех ядер из примеров 1-5, приведенных выше, при
любой такой последовательности и $\delta\in (0,\pi]$
$$
\lim_{n\to \infty}\int_0^\delta \sum _{s=0}^\infty (-1)^s\int_{x_{s+1}}^{x_s}K_n(t)dt=0.
$$

\bigskip


\begin{thebibliography}{99}

\bibitem{1} A. Kolmogoroff, {\it Une serie de Fourier-Lebesgue divergente
partout}, C.R. Acad. Sci. Paris, Ser.I {\bf 183} (1926), 1327--1329.

\bibitem{2} R. A. Hunt, {\it On the convergence of Fourier series}, Orthogonal Expansions and their Continuous Analogues
(Proc. Conf., Edwardsville, Ill., 1967), Southern Illinois Univ. Press, Carbondale, Ill, 1968, 235--255.

\bibitem{3} С.В. Конягин.  О расходимости всюду тригонометрических рядов
Фурье.  Матем. сб.,  191:1 (2000),  103--126 (Russian). - English transl. in  S.V. Konyagin, "On everywhere divergence of
trigonometric Fourier series",  Sb. Math. 191:1 (2000),  97--120.

\bibitem{4} N.Yu. Antonov, {\it Convergence of Fourier series},  East J. Appr., 2:2 (1996), 187-196.

\bibitem{5} A. Zygmund, {\it Trigonometric series}, v.I, II. 2nd ed. Cambridge
Univ. Press. N-Y, 1959, xii+383p.p., vii+354p.p.

\bibitem{6} G.H. Hardy and W.W. Rogosinski, {\it Fourier series},  2nd ed. Cambridge Tracts in Mathematics end
Mathematical Physics, 38,  Cambridge Univ. Press,  Cambridge,  1950,  x+100 p.p.

\bibitem{7} R.  Trigub and E.  Belinsky,  \emph{Fourier Analysis and Approximation of Functions},
Kluwer-Springer,  2004.

\bibitem{8} N.K. Bari, A Treatise in Trigonometric series,  Fizmatqiz, Moscow (Russian). - English transl. in Pergamon
Press, Mac. Millan, N-Y, 1964.

\bibitem{9} Д.Л. Фаддеев, {\it О представлении суммируемых функций сингулярными интегралами в точках
Lebesque'a}, Матем. сб., I (43) (1936), 351--368 (Russian).

\bibitem{10} G.H. Hardy, {\it Divergent series}, Oxford, 1949.

\bibitem{11} H. Lebesgue, {\it Recherches sur la convergence des series de
Fourier}, Math. Ann. Berlin-G\"ottingen-Heidelberg {\bf 61} (1905), 251--280.

\bibitem{12} Р.М. Тригуб, Суммируемость тригонометрических рядов Фурье в
$d$-точках, Изв. РАН, с.м., 79:4 (2015) (Russian).

\bibitem{13} E. Liflyand, S. Samko, and R. Trigub, {\it The Wiener algebra of absolutely convergent Fourier
integrals: an overview}, Anal. Math. Phys. {\bf 2} (2012), 1--68.

\bibitem{14} R.M. Trigub, {\it A Lower Bound for the $L_1$-norm of Fourier Series of Power
Type}, Mat. Zametki {\bf 73} (2003), 951--953 (Russian). - English transl. in Math. Notes {\it 73} (2003), 900--903.


\end{thebibliography}
\end{document}